\documentclass[11pt]{amsart}
\usepackage{amssymb,amsmath,amsfonts,amscd,euscript}

\newcommand{\nc}{\newcommand}

\nc{\slt}{\mathfrak{sl}_2}
\nc{\suth}{\widehat{\mathfrak{su}}(2)}
\nc{\g}{\mathfrak{g}}
\nc{\gh}{\widehat{\mathfrak{g}}}
\nc{\h}{\mathfrak{h}}
\nc{\la}{\lambda}
\nc{\slth}{\widehat{\slt}}
\nc{\C}{\mathbb C }
\nc{\Z}{\mathbb Z }
\nc{\N}{\mathbb N }
\nc{\al}{\alpha}
\nc{\be}{\beta}
\nc{\ve}{\varepsilon}
\nc{\ch}{{\mathop {\rm ch}}}
\nc{\Id}{{\mathop {\rm Id}}}
\nc{\Tr}{{\mathop {\rm Tr}\,}}
\nc{\U}{{\mathop {\rm U}}}
\nc{\bra}{\langle}
\nc{\ket}{\rangle}
\nc{\ld}{\ldots}
\nc{\cd}{\cdots}
\nc{\hk}{\hookrightarrow}
\nc{\n}{\mathfrak{n}}
\nc{\nh}{\widehat{\mathfrak{n}}}
\nc{\un}{\mathfrak{u}}
\nc{\T}{\otimes}
\nc{\bm}{{\bf m}}
\nc{\bs}{{\bf s}}
\nc{\bv}{{\bf v}}
\nc{\bt}{{\bf t}}
\nc{\bu}{{\bf u}}
\nc{\wt}{\widetilde}
\nc{\bin}[2]{{\genfrac{[}{]}{0pt}{0}{#1}{#2}}_q}
\nc{\fac}[1]{(#1)_q!}

\newtheorem{theo}{Theorem}[section]
\newtheorem{lem}{Lemma}[section]
\newtheorem{prop}{Proposition}[section]
\newtheorem{cor}{Corollary}[section]
\newtheorem{rem}{Remark}[section]

\begin{document}
\author{E.Feigin}
\title
[Infinite fusion products and $\slth$ cosets]
{Infinite fusion products and $\slth$ cosets}

\address{Evgeny Feigin:\newline
{\it Tamm Theory Division, Lebedev Physics Institute,
Russian Academy of Sciences,\newline
Russia, 119991, Moscow, Leninski prospect, 53}\newline
and \newline
{\it Independent University of Moscow,\newline
Russia, Moscow, 119002, Bol'shoi Vlas'evski, 11}}
\email{evgfeig@mccme.ru}

\begin{abstract}
In this paper we study an approximation of tensor product of
irreducible integrable $\widehat{\mathfrak{sl}_2}$ representations by
infinite fusion products. This gives an approximation of the
corresponding coset theories. As an application we represent characters of
spaces of these theories  as limits of certain restricted Kostka polynomials.
This leads
to the bosonic (which is known) and fermionic (which is new) formulas for the
$\widehat{\mathfrak{sl}_2}$ branching functions.
\end{abstract}
\maketitle

\section*{Introduction}
Let $\g$ be a semisimple Lie algebra, $\gh$ be the corresponding affine
algebra,
$$\gh=\g\T\C[t,t^{-1}]\oplus\C K\oplus \C d,$$
where $K$ is a central element and $[d,x_i]=-ix_i$. We set 
$\gh'=[\gh,\gh]=\g\T\C[t,t^{-1}]\oplus\C K$.
Let $L_{\la_1}$, $L_{\la_2}$ be two integrable irreducible $\gh$-modules.
Then one has the decomposition of the tensor product
$$L_{\la_1}\T L_{\la_2}=\bigoplus_\mu C_{\la_1,\la_2}^\mu \T L_\mu$$
into the direct sum of integrable irreducible representations of $\gh'$
(see \cite{Kac}). Here $C_{\la_1,\la_2}^\mu$
are spaces of highest weight vectors of the weight $\mu$ in the tensor
product $L_{\la_1}\T L_{\la_2}$. Therefore $C_{\la_1,\la_2}^\mu$ are naturally
graded by the operator $d$. This gives the character (branching function)
\begin{equation}
\label{branch}
c_{\la_1,\la_2}^\mu(q)=\Tr q^d|_{C_{\la_1,\la_2}^\mu}.
\end{equation}
Note that the GKO construction (see \cite{GKO}) endows
spaces $C_{\la_1,\la_2}^\mu$ with the structure of the representation of
the Virasoro algebra $Vir$. This also gives a grading which differs from
$(\ref{branch})$ by certain constant.

There exist different formulas for $c_{\la_1,\la_2}^\mu(q)$.
For the case $\g=\slt$ the bosonic (alternating sign) formula was obtained
in \cite{B, K, R} using
the representation theory of Virasoro algebra (Feigin-Fuchs construction 
\cite{FFu}). Another approach is
based on the connection of the branching functions with configuration sums
of RSOS model (see \cite{DJKMO}, \cite{S1}, \cite{S2}). 
This also gives different types formulas, in particular the fermionic formula 
for $\g=\slt$.
One of the important points in this approach is a construction of some
finitization (approximation) of branching functions. The same method is
used in \cite{SS}, where for the type $A$ affine Kac-Moody algebra the 
finitization of some branching functions is constructed by means of the 
combinatorics of
crystal bases. This allows to obtain $c_{\la_1,\la_2}^\mu(q)$ as certain
limits of restricted Kostka polynomials.  In our paper we construct the
representation theoretical approximation of the spaces $C_{\la_1,\la_2}^\mu$
for $\g=\slt$. We give some details below.

Let $L_{i,k}$, $0\le i\le k$ be irreducible integrable level $k$
representations with highest weight $i$ with respect to $h\T 1\in \slth$
($h$ is a standard generator of the Cartan subalgebra of $\slt$).
Then one has the isomorphism of $\slth'$ modules
$$L_{i_1,k_1}\T L_{i_2,k_2}=\bigoplus_{j=0}^{k_1+k_2} C_{i_1,i_2}^j\T
L_{j,k_1+k_2}.$$
Our main tool is a construction of the filtration
\begin{equation}
\label{filtr}
L(0)\hk L(1)\hk L(2)\hk\ldots =L_{i_1,k_1}\T L_{i_2,k_2}
\end{equation}
of the tensor product, where $L(p)$ are certain integrable representations of
$\slth$. Namely, let $v_n\in L_{i_1,k_1}$, $w_n\in L_{i_2,k_2}$, $n,m\in\Z$
be sets of extremal vectors (the orbits of the highest weight vectors
with respect to the action of the $\slth$ Weyl group). Then obviously
$$L_{i_1,k_1}\T L_{i_2,k_2}=\bigcup_{n,m\in\Z} U(\slth)\cdot (v_n\T w_m),$$
where $U(\slth)$ is the universal enveloping algebra. 
In addition it is easy to find $n(p), m(p)$, $p\ge 0$ such that for 
$L(p)=U(\slth)\cdot (v_{n(p)}\T w_{m(p)})$ the following holds:
\begin{equation}
\label{L_p}
L(0)\hk L(1)\hk L(2)\hk\ldots =L_{i_1,k_1}\T L_{i_2,k_2}.
\end{equation}
This procedure reduces the decomposition of the right hand side of
$(\ref{L_p})$ to the decomposition of $L(p)$ into the direct sum of 
irreducible representations. This can be done using the results from 
\cite{FF2, FF4}.

We recall that in \cite{FF2} the spaces $U(\slth)\cdot (v_n\T w_m)$ 
were
identified with infinite fusion products (the inductive limits of 
finite-dimensional fusion products, see \cite{FL}). 
The infinite fusion products were decomposed in \cite{FF4} and
the corresponding 
$q$-multiplicities of irreducible representations were expressed in terms 
of the restricted
Kostka polynomials. Therefore from
$(\ref{filtr})$ we obtain that branching functions are equal to the
appropriate limits of restricted Kostka polynomials.

This Kostka polynomial approximation gives two different formulas for branching
functions. From one hand we can use the alternating sum formula,
which expresses the restricted Kostka polynomials in terms of the unrestricted
Kostka polynomials (see \cite{SS, FJKLM}). 
The latter are related to the characters of the representations of $\slth$. 
Namely, 
certain limits of unrestricted Kostka polynomials can be expressed as a
difference of two $\slth$ string functions.
This leads to the following formula
\begin{multline*}
c_{i_1,i_2}^j(q)=q^{\gamma_1(i_1,i_2,j)}\times\\
\sum_{p\in\Z}
q^{-(k_1+k_2+2)p^2-(j+1)p} (\ch_q L_{i_1,k_1}^{2(k_1+k_2+2)p+j-i_2}-
\ch_q L_{i_1,k_1}^{2(k_1+k_2+2)p+j+i_2+2}),
\end{multline*}
where $L_{i,k}^a(q)=\{v\in L_{i,k}:\ (h\T 1)v= av\}$ 
and $\gamma_1(i_1,i_2,j)$ is some constant.
We show that this bosonic formula can be rewritten in a form of 
\cite{B, K, R}.

Another possibility is to use the fermionic formula for the restricted Kostka
polynomials (see \cite{SS, FJKLM}). In the appropriate limit this 
approach gives the following type formula:
\begin{multline}
c_{i_1,i_2}^j(q)= q^{\gamma_2(i_1,i_2,j)}\times\\
\sum_{\substack{s_i\ge 0 \\ i\in \{1,\ldots,k_1+k_2\}\setminus \{k_1\}}}
\frac{q^{\bs B \bs+\bu\bs}\bin{C\bs +\bv+\bs}{\bs}}
{\fac{\min(j,k_2)-i_2+ 2\sum_{\be\ne k_1} s_\be (\be-\min(k_1,\be))}},
\end{multline}
where $B$ and $C$ are some $(k_1+k_2-1)\times (k_1+k_2-1)$ matrixes and
$\bu,\bv$ are some vectors. The following notations are used: for two vectors
${\bf n}, {\bf m}\in \Z_{\ge 0}^N$ we set
$$\bin{\bf n}{\bf m}=\prod_{i=1}^N \bin{n_i}{m_i}=
\prod_{i=1}^N \frac{\fac{n_i}}{\fac{m_i}\fac{n_i-m_i}},\quad
\fac{k}=\prod_{i=1}^k (1-q^i).$$

Our paper is organized in the following way.

In Section $1$ we recall the main definitions and properties of the
representations of $\slth$ and of the fusion products.

Section $2$ is devoted to the description of the Kostka polynomial 
approach
to the computation of the branching functions.  This gives bosonic
(Section $3$) and fermionic (Section $4$) formulas for $\slth$ branching
functions.

{\bf Acknowledgements.} This work was partially supported by the RFBR grant
03-01-00167 and LSS 4401.2006.2.

\section{Preliminaries}
\subsection{Fusion product}
In this section we fix our notations and collect the main properties of
fusion products (see \cite{FL, FF1, FF2, FF3}).

Let $V_1,\ldots, V_n$ be cyclic representations of the Lie algebra $\g$ with
cyclic vectors $v_1,\ldots,v_n$. Fix
$Z=(z_1,\ldots,z_n)\in\C^n$ with $z_i\ne z_j$ for
$i\ne j$.  The fusion product
$V_1(z_1)*\ldots * V_n(z_n)$ is
the adjoint graded $\g\T \C[t]$ module with respect to the filtration $F_m$
on the tensor product $V_1(z_1)\T \cdots \T V_n(z_n)$:
\begin{equation}
\label{fusfiltr}
F_m=\mathrm{span}\{ g^{(1)}_{k_1}\cdots g^{(p)}_{k_p} (v_1\T\cdots\T v_n):\
k_1+\cdots +k_p\le m, g^{(i)}\in\g\}.
\end{equation}
Here $g_k=g\T t^k$ and $V_i(z_i)$ is the evaluation representation of
$\g\T\C[t]$, which is
isomorphic to $V_i$ as vector space and the action is defined via the
map $\g\T\C[t]\to\g$, $g\T t^j\mapsto z_i^j g$, $g\in\g$.

The most important property of fusion product is its independence on $Z$
in some special cases.
We will deal with the case $\g=\slt$.

Let $A=(a_1,\ldots,a_n)$. Denote by
\begin{equation}
\label{fus}
M_A=\pi_{a_1} *\ldots *\pi_{a_n}
\end{equation}
the fusion product of finite-dimensional
irreducible representations of $\slt$ ($\dim\pi_j=j+1$). Let $v_A$ be the
highest weight vector of $(\ref{fus})$ which is the image of the tensor
product of highest weight vectors of $\pi_{a_i}$.
We set
$$(\pi_{a_1} *\ldots *\pi_{a_n})_m=
\mathrm{span}\{(x^{(1)}_{i_1}\ldots x^{(k)}_{i_k})\cdot v_A,
i_1+\ldots +i_k=m, x^{(i)}\in\slt\}.$$
Let $h$ be the standard generator of the Cartan subalgebra of $\slt$. For
$\slt$-module $M$ we denote $M^\al=\{w\in M:\ hw=\al w\}$ and set
$$\ch_q  (\pi_{a_1} *\ldots *\pi_{a_n})^\al=
\sum_{m=0}^\infty q^m \dim
\bigl(
(\pi_{a_1} *\ldots *\pi_{a_n})^\al\cap (\pi_{a_1} *\ldots *\pi_{a_n})_m
\bigr).$$

We will need the following generalization of the standard embedding
$\pi_{a+b}\hk\pi_a\T\pi_b$. Let $A=(a_1\le\ldots\le a_n)$,
$B=(b_1\le\ldots\le b_m)$, $m\le n$. Then  we have an embedding
\begin{gather}
\label{comp}
\pi_{a_1}*\ldots * \pi_{a_{n-m}}*\pi_{a_{n-m+1}+b_1}*\ldots * \pi_{a_n+b_m}\hk
M_A\T M_B,\\
v_{a_1,\ldots, a_{n-m}, a_{n-m+1}+b_1, \ldots, a_n+b_m}\mapsto v_A\T v_B.
  \nonumber
\end{gather}

In the rest of this subsection we discuss one class of submodules of
fusion products (see \cite{FF3}).
Let $A=(a_1\le\ldots\le a_n)$.
Then for any $1\le i<  n$ there exists
$\slt\T\C[t]$-module $S_i(A)$ such that the following sequence is
exact:
\begin{multline}
0\to S_i(A)\to \pi_{a_1}*\ldots *\pi_{a_n}\to\\
\pi_{a_1}*\ldots * \pi_{a_{i-1}}* \pi_{a_i-1} *\pi_{a_{i+1}+1}* \pi_{a_{i+2}}
* \ldots * \pi_{a_n}\to 0
\end{multline}
For example, for $i=1$
$$S_1(A)\simeq \pi_{a_2-a_1}*\pi_{a_3}*\ldots *\pi_{a_n}.$$
We will also need the case $i=n-1$. In this case
\begin{equation}
\label{Sn-1}
S_{n-1}(A)\simeq \pi_{a_1}* \ldots *\pi_{a_{n-2}}\T \pi_{a_n-a_{n-1}}.
\end{equation}
Therefore one has an exact sequence
\begin{multline}
\label{lSn-1}
0\to \pi_{a_1}*\ldots *\pi_{a_{n-2}}\T \pi_{a_n-a_{n-1}}
\to \pi_{a_1}*\ldots *\pi_{a_n}\to \\
\pi_{a_1}*\ldots * \pi_{a_{n-2}}* \pi_{a_{n-1}-1}* \pi_{a_n+1}\to 0
\end{multline}

\subsection{Representations of $\slth$.}
In this subsection we fix our notations about $\slth$. Let
$$\slth=\slt\T\C[t,t^{-1}]\oplus\C K\oplus\C d,$$
where $K$ is a central element and $[d,x_i]=-ix_i$, where for
$x\in\slt$ we put $x_i=x\T t^i$.
Let $e,h,f$ be standard basis of $\slt$. Consider nilpotent subalgebras
$$\n_+=\slt\T t^{-1}\C[t^{-1}]\oplus\C f,\ 
  \n_-=\slt\T t\C[t]\oplus\C e.$$
We denote by $L_{l,k}$, $0\le l\le k$ integrable irreducible $\slth$-module
with highest weight vector $v_{l,k}$ such that
$$h_0v_{l,k}=lv_{l,k},\ K v_{l,k}=kv_{l,k},\ dv_{l,k}=0,\ \n_-v_{l,k}=0,\
U(\n_+)\cdot v_{l,k}=L_{l,k},$$ where
$U(\n_+)$ is the universal enveloping algebra.
Representations $L_{l,k}$ are bi-graded by operators $d$ and $h_0$.
We set
$$L_{l,k}=\bigoplus_{\al,s\in\Z} (L_{l,k})^\al_s=
\bigoplus_{\al,s\in\Z}\{v:\ dv=sv, h_0v=\al v\}.$$
This determines the character
$\ch_{q,z} L_{l,k}=\sum_{\al,s\in\Z} q^s z^\al \dim (L_{l,k})^\al_s.$
For any graded subspace $V\hk L_{l,k}$ we set
$$\ch_q V=\sum_{s\ge 0} q^s \dim \{v:\ dv=sv\}.$$

We now recall the Sugawara construction for the representation of the
Virasoro algebra in the space of the level $k$ $\slth$-module.
Namely, following operators form the Virasoro algebra:
$$L_n=\frac{1}{2(k+2)} \sum_{m\in\Z} :e_m f_{n-m} + f_m e_{n-m} +
\frac{1}{2} h_m h_{n-m}:, \ n\in\Z,$$
where $:\ :$ is the normal ordering
sign,
$$
:\ x_i y_j\ :=\begin{cases} x_i y_j, \ i\ge j,\\
y_j x_i, \ j\ge i,\\ \frac{1}{2} (x_i y_j+ y_j x_i).
\end{cases}
$$
The central charge is equal to $c(k)=\frac{3k}{k+2}$.
We denote by $\triangle_{l,k}$ the conformal weight of
the highest weight vector $v_{l,k}\in L_{l,k}$, i.e.
$L_0 v_{l,k}=\triangle_{l,k} v_{l,k}$.

We consider the
decomposition of the tensor product
$$L_{i_1,k_1}\T L_{i_2,k_2}=\bigoplus_{j=0}^{k_1+k_2}
C_{i_1,i_2}^j \T L_{j,k_1+k_2},$$
where $C_{i_1,i_2}^j$ is the space of highest weight vectors of the
weight $j$ in the tensor product 
$L_{i_1,k_1}\T L_{i_2,k_2}$.
Then by the GKO construction (see \cite{GKO}) the space  
$C_{i_1,i_2}^j$  is
a representation of the Virasoro algebra
$$L_n=L^{(1)}_n\T \Id + \Id\T L^{(2)}_n-L^{(diag)}_n,$$
where $L^{(1)}_n, L^{(2)}_n, L^{(diag)}_n$ are the Sugawara operators acting
on $L_{i_1,k_1}$, $L_{i_2,k_2}$ and  $L_{i_1,k_1}\T L_{i_2,k_2}$ respectively.
We put
$$c_{i_1,i_2}^j(q)=\ch_q C_{i_1,i_2}^j=
\Tr q^{L_0}|_{C_{\la_1,\la_2}^\mu}.$$
\begin{rem}
Note that the degree operator $d\in\gh$ acts on $C_{\la_1,\la_2}^\mu$ and
$L_0=d+\triangle_{i_1,k_1}+ \triangle_{i_2,k_2}-\triangle_{j,k_1+k_2}$.
\end{rem}

\subsection{Kostka polynomials and limit constructions.}
For the $k$-tuple
$\bm=(m_1,\ldots,m_k)\in\Z_{\ge 0}^k$  we set
\begin{equation*}
V_\bm=\pi_1^{*m_1}*\ldots *\pi_k^{*m_k}=
\underbrace{\pi_1*\ldots *\pi_1}_{m_1}*\ldots *
\underbrace{\pi_k*\ldots *\pi_k}_{m_k}.
\end{equation*}
We also use the notation $\bm=(1^{m_1}\ldots k^{m_k})$.

Consider the decomposition of $V_\bm$ into the direct sum of irreducible
representations of $\slt\T 1\hk\slt\T\C[t]$:
$$V_\bm=\bigoplus_{l\ge 0} \pi_l\T \wt K_{l,\bm},$$
where $\wt K_{l,\bm}\hk V_\bm$ is a space of highest weight vectors of weight 
$l$. We note that each $\wt K_{l,\bm}$ inherits a grading from $V_\bm$.
It is proved in \cite{FF4} that
$\ch_q \wt K_{l,\bm}=\wt K_{l,\bm}(q)$, where is unrestricted Kostka
polynomial. These polynomials are related to ones from \cite{FJKLM} by
$$\wt K_{l,\bm}(q)= q^{h(\bm)} K_{l,\bm}(q^{-1}),$$
where
\begin{gather}
\label{h}
h(\bm)=(\bm A\bm-p(\bm))/4, \quad A= (A_{i,j})=\min(i,j), \\ \nonumber
p(\bm)=\#\{\al:\ m_\al+\ldots +m_k \text{ is odd }\}.
\end{gather}
We note that $h(\bm)$ can be defined as follows:
$$h(\bm)=\max\{s:\ (\pi_1^{*m_1}*\ldots *\pi_k^{*m_k})_s\ne 0\}.$$
Thus the "reversed" character of the fusion product is given by
$$\wt\ch_q V^\al_\bm=q^{h(\bm)} \ch_{q^{-1}} V^\al_{\bm}.$$

We proceed with a limit construction of fusion products.
It is proved in \cite{FF2} that there exists an embedding of
$\slt\T\C[t]$-modules:
$$V_{(m_1,\ldots,m_k)}\hk V_{(m_1,\ldots, m_{k-1}, m_k+2)}.$$
This allows to define
an injective limit
$$L_{\bm,k}=\lim_{N\to\infty} V_{(m_1,\ldots, m_{k-1}, m_k+2N)}.$$
It turns out that $L_{\bm,k}$  has the natural structure of level $k$ 
$\slth$ module.
For example, $L_{i,k}=\lim_{N\to\infty} \pi_i*\pi_k^{*2N}$.
In general representations $L_{\bm,k}$ are reducible. Consider the
decomposition
$$L_{\bm,k}=\bigoplus_{l=0}^k C_{l,\bm}\T L_{l,k},$$
where $C_{l,\bm}\hk L_{\bm,k}$ is the space of the highest weight vectors of
the weight $l$.  Note that $C_{l,\bm}$ are naturally garded by the operator
$d$. We set $\ch_q C_{l,\bm}=\Tr q^{d}|_{C_{l,\bm}}$. It is shown in
\cite{FF4} that
\begin{equation}
\label{decomp}
\ch_q C_{l,\bm}=\wt K_{l,\bm}^{(k)}(q),
\end{equation}
where $\wt K^{(k)}_{l,\bm}(q)$ is the restricted Kostka polynomials.
These $\wt K_{l,\bm}^{(k)}(q)$ are relates to
  $K_{l,\bm}^{(k)}(q)$  from \cite{FJKLM} by
$$\wt K_{l,\bm}^{(k)}(q)= q^{h(\bm)} K_{l,\bm}^{(k)}(q^{-1}).$$

\section{Kostka polynomials approximation}
In this section we obtain the Kostka polynomials approximation of the
characters of the $\slth$ coset  models, using the injective limits of
fusion products.

Let $v_n\in L_{i_1,k_1}$, $n\in\Z$ be the set of extremal vectors,
$v_n=T^n v_{i,k}$, where $T$ is the translation operator from the Weyl group of
$\slth$ and $v_{i,k}$ is a highest weight vector. We have
$h_0 v_n=i_1-2k_1 n$. We also denote the set of extremal vectors of
$L_{i_2,k_2}$ by $w_n=T^n v_{i_2,k_2}$, $h_0 w_n=i_2-2k_2n$.

\begin{lem}
\label{1}
Let $n\ge m$. Then\\
$a)$. \ $U(\slth)\cdot (v_n\T w_m)\hk U(\slth)\cdot (v_{n+1}\T w_m)$;\\
$b)$. \ $\lim_{n\to\infty} U(\slth)\cdot (v_n\T w_m)=L_{i_1,k_1}\T L_{i_2,k_2}.$
\end{lem}
\begin{proof}
We note that
$(e_{2n-1})^{k_1-i_1} (e_{2n})^{i_1} v_n$ is proportional to $v_{n-1}$ and
also
$e_i w_m=0$ if $i> 2m$. Therefore
$(e_{2n-1})^{k_1-i_1} (e_{2n})^{i_1} (v_n\T w_m)$ is proportional to
$v_{n-1}\T w_m$. So $a)$ is proved.

To prove $b)$ it suffices to note that
$$U(\slth)\cdot \mathrm{span} \{ v_n\T w_m:\ n\ge m\}= L_{i_1,k_1}\T
L_{i_2,k_2}.$$
\end{proof}

\begin{lem}
\label{2}
Let $n\ge m> 0$. Then
\begin{equation}
\label{lim}
U(\slth)\cdot (v_n\T w_m)\simeq \lim_{N\to\infty}
\pi_{i_1} * \pi_{k_1}^{*(2(n-m)-1)} *\pi_{k_1+i_2} *\pi_{k_1+k_2}^{*2N}.
\end{equation}
\end{lem}
\begin{proof}
We note that
$$U(\slt\T\C[t])\cdot v_n\simeq \pi_{i_1} * \pi_{k_1}^{*2n},\
   U(\slt\T\C[t])\cdot w_m\simeq \pi_{i_2} * \pi_{k_2}^{*2m}.$$
Therefore, because of $(\ref{comp})$,
$$U(\slt\T\C[t])\cdot (v_n\T w_m)\simeq
\pi_{i_1} * \pi_{k_1}^{* (2(n-m)-1)} *\pi_{k_1+i_2} *\pi_{k_1+k_2}^{*2m}.$$
We now obtain our lemma from
$$U(\slth)\cdot (v_n\T w_m)\simeq
\lim_{N\to\infty} U(\slt\T\C[t])\cdot (v_{n+N}\T w_{m+N}).$$
\end{proof}

We denote by $\bm(N)$ a $(k_1+k_2)$-tuple such that
\begin{equation}
\label{m(N)}
V_{\bm(N)}=
\pi_{i_1} * \pi_{k_1}^{* (2N-1)} *\pi_{k_1+i_2}.
\end{equation}
From $(\ref{lim})$  and Lemma $\ref{1}$ we obtain
\begin{equation}
\label{appr}
L_{i_1,k_1}\T L_{i_2,k_2}=\lim_{N\to\infty}
L_{\bm(N), k_1+k_2}.
\end{equation}
Now $(\ref{decomp})$ and $(\ref{appr})$ gives
\begin{cor}
$$c_{i_1,i_2}^j(q)=
q^{\triangle_{i_1,k_1}+\triangle_{i_2,k_2}-\triangle_{j,k_1+k_2}}
\lim_{N\to\infty} \wt K^{(k_1+k_2)}_{j,\bm(N)}.$$
\end{cor}

\section{Bosonic formula}
We use the alternating sign formula for the restricted Kostka polynomials
in terms of the unrestricted Kostka polynomials (see \cite{SS, FJKLM})
\begin{multline}
K^{(k)}_{j,\bm}(q)=\\
\sum_{p\ge 0} q^{(k+2)p^2+(j+1)p} K_{2(k+2)p+j,\bm}(q)-
                    \sum_{p> 0} q^{(k+2)p^2-(j+1)p} K_{2(k+2)p-j-2,\bm}(q).
\end{multline}
Recall the notations from \cite{FF4}
$$\wt K^{(k)}_{l,\bm}(q)=q^{h(\bm)} K^{(k)}_{l,\bm}(q^{-1}),\
   \wt K_{l,\bm}(q)=q^{h(\bm)} K_{l,\bm}(q^{-1}).$$

\begin{lem}
\begin{equation}
\label{conv}
\lim_{N\to\infty} \wt K_{j,\bm(N)}(q)=\ch_q L_{i_1,k_1}^{j-i_2}-
                            \ch_q L_{i_1,k_1}^{j+i_2+2}.
\end{equation}
\end{lem}
\begin{proof}
It is shown in \cite{FF4} that $\wt K_{j,\bm}(q)$ is a multiplicity of 
$\pi_j$ in $V_\bm$.
Consider embeddings
\begin{equation}
\label{mil}
\pi_{i_1}* \pi_{k_1}^{*2(N-1)}\T \pi_{i_2} \hk \pi_{i_1} * \pi_{k_1}^{*(2N-1)}
* \pi_{k_1+i_2}\hk  \pi_{i_1}* \pi_{k_1}^{*2N}\T \pi_{i_2},
\end{equation}
where the first embedding comes from $(\ref{comp})$ and the second from
$(\ref{Sn-1})$.
We note that $(\ref{mil})$ for $N$ and $N+1$ can be combined into the
commutative diagram
$$
\begin{CD}
\pi_{i_1}* \pi_{k_1}^{*2(N-1)}\T \pi_{i_2} @> > > 
\pi_{i_1} * \pi_{k_1}^{*(2N-1)}
* \pi_{k_1+i_2}@> > >   \pi_{i_1}* \pi_{k_1}^{*2N}\T \pi_{i_2}\\
@VVV    @.   @VVV \\
\pi_{i_1}* \pi_{k_1}^{*2N}\T \pi_{i_2} @> > >  \pi_{i_1} * \pi_{k_1}^{*(2N+1)}
* \pi_{k_1+i_2} @> > >  \pi_{i_1}* \pi_{k_1}^{*2(N+1)}\T \pi_{i_2}.
\end{CD}
$$

In view of $L_{i,k}=\lim_{N\to\infty} \pi_i* \pi_k^{*2N}$ we obtain
\begin{multline}
\wt K_{j,\bm(N)}(q)=\wt\ch_q (V_{\bm(N)}\T \pi_{i_2})^j-
                     \wt\ch_q (V_{\bm(N)}\T \pi_{i_2})^{j+2} =\\
\wt\ch_q V_{\bm(N)}^{j-i_2}-
                     \wt\ch_q V_{\bm(N)}^{j+i_2+2}
\to
\ch_q L_{i_1,k_1}^{j-i_2}- \ch_q L_{i_1,k_1}^{j+i_2+2},\ N\to\infty.
\end{multline}
Lemma is proved.
\end{proof}

Now we need to know how "fast" the left hand side of $(\ref{conv})$
converges to the right hand side.

\begin{prop}
\label{fast}
$\ch_q L_{i,k}^n -\wt\ch_q (\pi_i * \pi_k^{*2N})^n=
\mathrm{O}(q^{N+1+\frac{n^2}{4k}-\frac{k}{4}})$.
\end{prop}
\begin{proof}
We first consider the difference
$$\wt\ch_q (\pi_i * \pi_k^{*2(N+1)})^n -\wt\ch_q (\pi_i * \pi_k^{*2N})^n.$$
Recall (see $(\ref{lSn-1})$) that there exists an exact sequence
$$0\to  \pi_i * \pi_k^{*2N}\to \pi_i * \pi_k^{*2(N+1)}\to
\pi_i* \pi_{k-1} *\pi_k^{*2N} *\pi_{k+1}\to 0.$$
We also note that the $\slt\T\C[t]$-homomorphism
$$\pi_i * \pi_k^{*2(N+1)}\to\pi_i* \pi_{k-1} *\pi_k^{*2N} *\pi_{k+1}$$
is determined by the condition that the highest weight vector (with respect
to $h_0$) maps to the
highest weight vector. Therefore
\begin{multline}
\wt\ch_q (\pi_i * \pi_k^{*2(N+1)})^n -\wt\ch_q (\pi_i * \pi_k^{*2N})^n= \\
q^{h(i^1 k^{2(N+1)})-h(i^1 (k-1) k^{2N} (k+1))}
\wt\ch_q (\pi_i* \pi_{k-1} *\pi_k^{*2N} *\pi_{k+1})^n.
\end{multline}
Using the formula $h(\bm)=(\bm A\bm-p(\bm))/4$ (see $(\ref{h})$) we obtain
\begin{multline}
\label{hdif}
h(i^1 k^{2(N+1)})-h(i^1 (k-1) k^{2N} (k+1))=\\
(N+1)(k(N+1)+i)-(N+1)(k(N+1)+i-1)=N+1.
\end{multline}
To evaluate $\wt\ch_q (\pi_i* \pi_{k-1} *\pi_k^{*2N} *\pi_{k+1})^n$ we use
the following formula for the character of the fusion product from \cite{FF1}
\begin{multline}
\wt\ch_q (\pi_1^{*m_1}*\ldots *\pi_k^{*m_k})^n=\\
q^{-\frac{p(\bm)}{4}}
\sum_{\substack{j_1,\ldots,j_k\ge 0\\
2\sum_{l=1}^k (j_l-\al_l)=n}}
q^{\sum\nolimits_{l=1}^k (j_l-\al_l)^2} \bin{m_k}{j_k}\prod_{l=1}^{k-1}
\bin{m_{k-l}+j_{k-l+1}}{j_{k-l}},
\end{multline}
where $2\al_l=m_l+\ldots + m_k$, and 
$$\bin{m}{j}=\frac{\fac{m}}{\fac{j}\fac{m-j}},\ 
\fac{m}=\prod_{i=1}^m (1-q^i).$$
In view of
$\sum_{l=1}^k (j_l-\al_l)=n/2$ we obtain
$$-\frac{p(\bm)}{4}+ \sum_{l=1}^k (j_l-\al_l)^2\ge \frac{n^2}{4k}-
\frac{k+1}{4}.$$
Therefore
\begin{equation}
\label{ndif}
\wt\ch_q (\pi_i* \pi_{k-1} *\pi_k^{*2N} *\pi_{k+1})^n=
\mathrm{O}(q^{\frac{n^2}{4k}-\frac{k}{4}}).
\end{equation}
From $(\ref{hdif})$ and $(\ref{ndif})$ we obtain
\begin{equation}
\wt\ch_q (\pi_i * \pi_k^{*2(N+1)})^n -\wt\ch_q (\pi_i * \pi_k^{*2N})^n=
\mathrm{O}(q^{N+1+\frac{n^2}{4k}-\frac{k}{4}}).
\end{equation}
Now using the limit construction
$$L_{i,k}=\pi_i\hk \ldots\hk \pi_i*\pi_k^{*2N}\hk \pi_i*\pi_k^{*2(N+1)}\hk
\ldots$$
we obtain our proposition.
\end{proof}

\begin{cor}
Let $l\ge i_2$.
\label{O}
$$\wt K_{l,\bm(N)}(q)-(\ch_q L_{i_1,k_1}^{l-i_2}- \ch_q L_{i_1,k_1}^{l+i_2+2})=
\mathrm{O}(q^{N+\frac{(l-i_2)^2}{4k_1}-\frac{k_1}{4}}).$$
\end{cor}
\begin{proof}
We recall that
\begin{equation}
\label{K}
\wt K_{l,\bm(N)}(q)=\wt\ch_q (\pi_{i_1}*\pi_{k_1}^{*(2N-1)}*\pi_{k_1+i_2})^l-
\wt\ch_q (\pi_{i_1}*\pi_{k_1}^{*(2N-1)}*\pi_{k_1+i_2})^{l+2}.
\end{equation}
Using $(\ref{comp})$ and $(\ref{Sn-1})$ we obtain embeddings
\begin{equation}
\label{apr}
\pi_{i_1}*\pi_{k_1}^{*(2N-2)}\T \pi_{i_2}\hk
\pi_{i_1}*\pi_{k_1}^{*(2N-1)}*
\pi_{k_1+i_2}\hk \pi_{i_1}*\pi_{k_1}^{*(2N)} \T \pi_{i_2}.
\end{equation}
From Proposition $(\ref{fast})$ we have
\begin{multline}
\label{sum}
\wt \ch_q (\pi_{i_1}*\pi_{k_1}^{*(2N)} \T \pi_{i_2})^l-
\ch_q (L_{i_1,k_1}\T \pi_{i_2})^l=\\
\sum_{s=0}^{i_2}
\bigl(\wt \ch_q (\pi_{i_1}*\pi_{k_1}^{*(2N)} \T \pi_{i_2})^{l-i_2+2s}-
\ch_q (L_{i_1,k_1}\T \pi_{i_2})^{l-i_2+2s}\bigr)=\\
\mathrm{O} (q^{N+1+\frac{(l-i_2)^2}{4k_1}-\frac{k_1}{4}}),
\end{multline}
because $l\ge i_2$.
Therefore from $(\ref{apr})$ and $(\ref{sum})$ follows that
\begin{equation}
\wt \ch_q (\pi_{i_1}*\pi_{k_1}^{*(2N-1)}*\pi_{k_1+i_2})^l-
\ch_q (L_{i_1,k_1}\T \pi_{i_2})^l=
\mathrm{O} (q^{N+\frac{(l-i_2)^2}{4k_1}-\frac{k_1}{4}}),
\end{equation}
where the approximation $L_{i_1,k_1}\T\pi_{i_2}=
\lim_{N\to\infty} \pi_{i_1}* \pi_{k_1}^{*2N} \T\pi_{i_2}$ is used.
Now our corollary follows from  $(\ref{K})$ and
$$\ch_q (L_{i_1,k_1}\T\pi_{i_2})^l-\ch_q (L_{i_1,k_1}\T\pi_{i_2})^{l+2}=
\ch_q L_{i_1,k_1}^{l-i_2}- \ch_q L_{i_1,k_1}^{l+i_2+2}.$$
\end{proof}

\begin{theo}
\begin{multline}
\label{bos}
\lim_{N\to\infty} \wt K^{(k_1+k_2)}_{j,\bm(N)}=\\ \sum_{p\in\Z}
q^{-(k_1+k_2+2)p^2-(j+1)p}\bigl( \ch_q L_{i_1,k_1}^{2(k_1+k_2+2)p+j-i_2}-
\ch_q L_{i_1,k_1}^{2(k_1+k_2+2)p+j+i_2+2}\bigr).
\end{multline}
\end{theo}

\begin{proof}
We use the alternating sign formula
\begin{multline}
\wt K^{(k_1+k_2)}_{j,\bm(N)}(q)=
\sum_{p\ge 0} q^{-(k_1+k_2+2)p^2-(j+1)p} \wt K_{2(k_1+k_2+2)p+j,\bm(N)}(q)-\\
    \sum_{p> 0} q^{-(k_1+k_2+2)p^2+(j+1)p} \wt K_{2(k_1+k_2+2)p-j-2,\bm(N)}(q).
\end{multline}
Using Corollary $(\ref{O})$ we rewrite this expression as
\begin{multline}
\sum_{p\ge 0} q^{-(k_1+k_2+2)p^2-(j+1)p}
\Bigl(\ch_q L_{i_1,k_1}^{2(k_1+k_2+2)p+j-i_2}- \\
  \ch_q L_{i_1,k_1}^{2(k_1+k_2+2)p+j+i_2+2}+
\mathrm{O}\bigl(q^{N+\frac{(2(k_1+k_2+2)p+j-i_2)^2-k_1^2}{4k_1}}\bigr)
\Bigr)-\\
\shoveleft {\sum_{p>  0} q^{-(k_1+k_2+2)p^2+(j+1)p}
\Bigl(\ch_q L_{i_1,k_1}^{2(k_1+k_2+2)p-j-2-i_2}-}\\
  \ch_q L_{i_1,k_1}^{2(k_1+k_2+2)p-j+i_2}+
\mathrm{O}\bigl(q^{N+\frac{(2(k_1+k_2+2)p-j-2-i_2)^2-k_1^2}{4k_1}}
\bigr)\Bigr).
\end{multline}
We note that for $p$ big enough we have
\begin{gather*}
(k_1+k_2+2)p^2+(j+1)p< \frac{(2(k_1+k_2+2)p-j-2-i_2)^2-k_1^2}{4k_1},\\
(k_1+k_2+2)p^2-(j+1)p<  \frac{(2(k_1+k_2+2)p-j-2-i_2)^2-k_1^2}{4k_1}.
\end{gather*}
Therefore we obtain
\begin{multline}
\lim_{N\to\infty} \wt K^{(k_1+k_2)}_{j,\bm(N)}(q)= \\
\sum_{p\ge 0} q^{-(k_1+k_2+2)p^2-(j+1)p}
\Bigl(\ch_q L_{i_1,k_1}^{2(k_1+k_2+2)p+j-i_2}-
  \ch_q L_{i_1,k_1}^{2(k_1+k_2+2)p+j+i_2+2}\Bigr)-\\
\sum_{p>  0} q^{-(k_1+k_2+2)p^2+(j+1)p}
\Bigl(\ch_q L_{i_1,k_1}^{2(k_1+k_2+2)p-j-2-i_2}-
  \ch_q L_{i_1,k_1}^{2(k_1+k_2+2)p-j+i_2}\Bigr).
\end{multline}
We now rewrite the second sum replacing  $p$ by $-p$.
This gives
\begin{equation}
-\sum_{p<  0} q^{-(k_1+k_2+2)p^2-(j+1)p}
\Bigl(\ch_q L_{i_1,k_1}^{2(k_1+k_2+2)p+j-i_2}-
  \ch_q L_{i_1,k_1}^{2(k_1+k_2+2)p+j+i_2+2}\Bigr)
\end{equation}
(because $\ch_q L_{i,k}^a=\ch_q L_{i,k}^{-a}$). Our theorem is proved.
\end{proof}

\begin{cor}
The right hand side of $(\ref{bos})$ coincides with
\begin{equation}
\label{coset}
q^{-\triangle_{i_1,k_1}- \triangle_{i_2,k_2}+ \triangle_{j,k_1+k_2}} 
c_{i_1,i_2}^j(q).
\end{equation}
\end{cor}

We finish this section with the identification of our bosonic formula
with the known one (see \cite{B},\cite{K},\cite{R}). 

Let $\h$ be Cartan subalgebra of $\slth$
and define elements $(a,k,s)\in \h^*$ by 
$$(a,k,s)h_0=i,\ (i,k,s)K=k,\ (a,k,s)d=s.$$
We consider a translation element $t$ from the Weyl group $W$ of $\slth$ 
defined by
$$t(a,k,s)=(a+2k,k,s+k+a).$$
Therefore we have an isomorphism of vector spaces 
$L_{i,k}^a\simeq L_{i,k}^{a+2k}$ and for the corresponding characters we obtain
\begin{equation}
\ch_q L_{i,k}^{a+2k}=q^{k+a} \ch_q L_{i,k}^a.
\end{equation}
This gives 
\begin{equation}
\label{shift}
\ch_q L_{i,k}^{a+2\la k}=q^{\la(\la k+a)} \ch_q L_{i,k}^a.
\end{equation}

We now rewrite the right hand side of $(\ref{bos})$ using $(\ref{shift})$.
Namely, let 
\begin{equation}
\label{change}
2(k_1+k_2+2)p+j-i_2=2m+2\la k_1,
\end{equation}
where $0\le m\le k_1/2$ is integer for even $i_1$ and half-integer for odd 
$i_1$. Using $(\ref{shift})$ and $(\ref{change})$  we obtain  
\begin{multline*}
q^{-(k_1+k_2+2)p^2-(j+1)p} \ch_q L_{i_1,k_1}^{2(k_1+k_2+2)p+j-i_2}=\\
q^{\frac{(j-i_2)^2}{4k_1^2}}q^{\frac{-m^2}{k_1}}
q^{\frac{p}{k_1}(p(k_1+k_2+2)(k_2+2)+(k_2+2)(j+1)-(k_1+k_2+2)(i_2+1))}
\ch_q L_{i_1,k_1}^m.
\end{multline*}
Combining this with the similar formula for the second term in the right hand 
side of $(\ref{bos})$ we obtain that up to a power of $q$ the
branching function $c_{i_1,i_2}^j(q)$ equals to
\begin{multline*}
\sum_{0\le m\le k_1/2} q^{-\frac{m^2}{4}} \ch_q L_{i_1,k_1}^{2m} \times\\
\Bigl( \sum_{\substack{p\in\Z \\ m_{j-i_2}(p)\equiv \pm m\pmod{k_1}}} 
q^{\frac{p}{k_1}(p(k_1+k_2+2)(k_2+2)+(k_2+2)(j+1)-(k_1+k_2+2)(i_2+1))} - \\
\sum_{\substack{p\in\Z \\ m_{j+i_2}(p)\equiv \pm m\pmod{k_1}}} 
q^{\frac{1}{k_1}(p(k_1+k_2+2)+j+1)(p(k_2+2)+(i_2+1))}\Bigr),
\end{multline*}
where $m_a(p)=a/2+(k_1+k_2+2)p$ and $m$ runs over integers if $i_1$ is even
and over half-integers if $i_1$ is odd.
Identifying $q^{-\frac{m^2}{4}} \ch_q L_{i_1,k_1}^{2m}$ with the $\slth$
string functions we obtain the known bosonic formula.

\section{Fermionic formula}
We now compute the limit $\lim_{N\to\infty} \wt K^{(k_1+k_2)}_{j,\bm(N)}(q)$
using the fermionic formula from \cite{FF4}:
\begin{equation}
\label{ferm}
q^{p(\bm)/4}\wt K^{(k)}_{j,\bm}=\sum_{\substack{\bs\in\Z_{\ge 0}^k\\ 2|\bs|=|\bm|-j}}
q^{(\frac{\bm}{2}-\bs)A (\frac{\bm}{2}-\bs)}\bin{A(\bm-2\bs)-\nu+\bs}{\bs},
\end{equation}
where $A=(A_{i,j})_{i,j=1}^k=(\min(i,j))$, $\nu=(\nu_a)_{a=1}^k=
(\max(0,a-k+j))$, $|\bm|=\sum_{i=1}^k im_i$. For $\bv\in \Z^k_{\ge 0}$ we
put
$$\bin{\bv}{\bs}=\prod_{i=1}^k \bin{v_i}{s_i}=\prod_{i=1}^k
\frac{\fac{v_i}}{\fac{s_i}\fac{v_i-s_i}}.$$
Now let $k=k_1+k_2$, $\bm=\bm(N)$ (see $(\ref{m(N)})$). Then
$|\bm(N)|=i_1+i_2+2k_1N$ and
for $\bs$ from the right hand side of $(\ref{ferm})$ we have
$|\bs|=k_1N+\frac{i_1+i_2-j}{2}$.
This gives
\begin{equation}
\label{sk1}
s_{k_1}=N+\frac{i_1+i_2-j}{2k_1}-\frac{1}{k_1}
\sum_{\substack{1\le\al\le k_1+k_2\\ \al\ne k_1}} \al s_\al.
\end{equation}
We now rewrite the fermionic formula for $\wt K^{(k_1+k_2)}_{j,\bm(N)}(q)$
using $(\ref{sk1})$.

We start with the power $(\frac{\bm(N)}{2}-\bs)A (\frac{\bm(N)}{2}-\bs)$.
Note that
$$\frac{\bm(N)}{2} A \frac{\bm(N)}{2}=N^2 k_1+ N i_1 +\frac{i_1+i_2}{4}.$$
Therefore
\begin{multline}
\label{A}
\left(\frac{\bm(N)}{2}-\bs\right)A \left(\frac{\bm(N)}{2}-\bs\right)=
N^2 k_1+ N i_1 +\frac{i_1+i_2}{4}+\\
\sum_{\substack{1\le \al,\be \le k_1+k_2\\ \al,\be\ne k_1}}
\min(\al,\be)s_\al s_\be - \\
\sum_{\al\ne k_1} s_\al (\min(\al,i_1)+(2N-1)\min(\al,k_1)+\min(\al, k_1+i_2))+
\\
k_1 s_{k_1}^2+ 2\sum_{\al\ne k_1}\min(\al,k_1) s_\al s_{k_1}-
(i_1+2Nk_1)s_{k_1}.
\end{multline}
Using $(\ref{sk1})$ we rewrite the last line as
\begin{multline}
\label{part}
k_1 (N+\frac{i_1+i_2-j}{2k_1}-\frac{1}{k_1}
\sum_{\al\ne k_1} \al s_\al)^2+ \\
(-i_1-2Nk_1+2\sum_{\al\ne k_1}\min(\al,k_1) s_\al)
(N+\frac{i_1+i_2-j}{2k_1}-\frac{1}{k_1}
\sum_{\al\ne k_1} \al s_\al).
\end{multline}
Combining together $(\ref{A})$ and $(\ref{part})$ we obtain
\begin{multline}
\left(\frac{\bm(N)}{2}-\bs\right)A \left(\frac{\bm(N)}{2}-\bs\right)=\\
\sum_{\al,\be\ne k_1} s_\al s_\be \left[ \min(\al,\be)+
\frac{1}{k_1}(\al\be-\min(\al,k_1)\be-\min(\be,k_1)\al)\right] +\\
\sum_{\al\ne k_1} s_\al \Bigl[-\min(\al,i_1)-\min(\al,k_1+i_2)+
\min(\al,k_1)\frac{i_1+i_2-j+k_1}{k_1}+\\ \al\frac{j-i_2}{k_1}\Bigr] +
\frac{i_1+i_2}{4}+\frac{(i_1+i_2-j)^2}{4k_1}-\frac{i_1(i_1+i_2-j)}{2k_1}.
\end{multline}

We now rewrite the binomial coefficient
\begin{equation}
\bin{(A(\bm(N)-2\bs))_\al-\nu_\al+s_\al}{s_\al}
\end{equation}
using $(\ref{sk1})$. Let $\al\ne k_1$. Then
\begin{multline}
(A(\bm(N)-2\bs))_\al-\nu_\al+s_\al=\\
\min(\al,i_1)+(2N-1)\min(\al,k_1)+\min(\al,k_1+i_2)-\\
-2\sum_{\be\ne k_1}\min(\al,\be) s_\be-2\min(\al,k_1)s_{k_1}-\nu_\al+s_\al
=\\
\min(\al,i_1)-\min(\al,k_1)(\frac{i_1+i_2-j+k_1}{k_1})+\min(\al,k_1+i_2)+\\
2\sum_{\be\ne k_1} s_\be \left(-\min(\al,\be)+\frac{\min(\al,k_1)\be}{k_1}
\right)- \nu_\al+s_\al.
\end{multline}
We note that the result is independent on $N$. Now let $\al=k_1$. Then
\begin{multline}
(A(\bm(N)-2\bs))_{k_1}-\nu_{k_1}=\\
i_1+(2N-1)k_1+k_1-2\sum_{\be\ne k_1} \min(k_1,\be)s_\be -2k_1s_{k_1}-
\nu_{k_1}=\\
j-i_2-\nu_{k_1}+ 2\sum_{\be\ne k_1} s_\be (\be-\min(k_1,\be)).
\end{multline}
Therefore
\begin{multline}
\bin{(A(\bm(N)-2\bs))_{k_1}-\nu_{k_1}+s_{k_1}}{s_{k_1}}=\\
\frac{1+O(q^{1+s_{k_1}})}
{\fac{j-i_2-\nu_{k_1}+ 2\sum_{\be\ne k_1} s_\be (\be-\min(k_1,\be))}}.
\end{multline}
Note that
\begin{multline}
s_{k_1}+\left(\frac{\bm(N)}{2}-\bs\right)A \left(\frac{\bm(N)}{2}-\bs\right)=\\
s_{k_1}+\sum_{\al=1}^{k_1+k_2}
\left[\sum_{\be=k_1+k_2-\al+1}^{k_1+k_2}
\left(\frac{\bm(N)}{2}-\bs\right)_\be\right]^2\ge \\
s_{k_1}+
\left[\sum_{\be=k_1+1}^{k_1+k_2} \left(\frac{\bm(N)}{2}-\bs\right)_\be\right]^2
+\left[\sum_{\be=k_1}^{k_1+k_2} \left(\frac{\bm(N)}{2}-\bs\right)_\be\right]^2=
\\
s_{k_1}+(s_{k_1+1}+\ldots +s_{k_1+k_2}-1/2)^2 +
(s_{k_1}+\ldots +s_{k_1+k_2}-N)^2.
\end{multline}
The expression in the last line is greater than or equal to $N/3$
(because if $s_{k_1}< N/3$ and $s_{k_1+1}+\ldots +s_{k_1+k_2}-1/2< N/3$ then
$(s_{k_1}+\ldots +s_{k_1+k_2}-N)^2> N/3$). Therefore
\begin{multline}
q^{(\frac{\bm}{2}-\bs)A (\frac{\bm}{2}-\bs)}
\bin{A(\bm-2\bs)_{k_1}-\nu_{k_1}+s_{k_1}}{s_{k_1}}=\\
\frac{1} {\fac{j-i_2-\nu_{k_1}+ 2\sum_{\be\ne k_1} s_\be (\be-\min(k_1,\be))}}
+ O(q^{N/3}).
\end{multline}

Using $p(\bm(N))=i_1+i_2$ we obtain
the following theorem.
\begin{theo}
\begin{multline}
\label{mf}
\lim_{N\to\infty}
q^{-\frac{(i_1+i_2-j)(i_2-i_1-j)}{4k_1}}
\wt K^{(k_1+k_2)}_{j,\bm(N)}= \\
\sum_{\substack{s_i\ge 0 \\ i\in \{1,\ldots,k_1+k_2\}\setminus \{k_1\}}}
\frac{q^{\bs B \bs+\bu\bs}\bin{C\bs +\bv+\bs}{\bs}}
{\fac{\min(j,k_2)-i_2+ 2\sum_{\be\ne k_1} s_\be (\be-\min(k_1,\be))}},
\end{multline}
where
\begin{multline*}
B_{\al,\be}= \min(\al,\be)+
\frac{\al\be-\min(\al,k_1)\be-\min(\be,k_1)\al}{k_1},\\
\shoveleft {\bu_\al=-\min(\al,i_1)-\min(\al,k_1+i_2)+} \\
\shoveright {\frac{\min(\al,k_1)(i_1+i_2-j+k_1)+ \al(j-i_2)}{k_1},}\\
\shoveleft{C_{\al,\be}=2\frac{\min(\al,k_1)\be}{k_1}-2\min(\al,\be)},\\
\shoveleft{\bv_\al=\min(\al,i_1)-\min(\al,k_1)\frac{i_1+i_2-j+k_1}{k_1}+}\\
\min(\al,k_1+i_2)-\max(0,\al-k_1-k_2+j).
\end{multline*}
\end{theo}

\begin{cor}
The branching function $c_{i_1,i_2}^j(q)$ equals to the product of the 
right hand side
of $(\ref{mf})$ and 
$$q^{\triangle_{i_1,k_1}+ \triangle_{i_2,k_2}-\triangle_{j,k_1+k_2}+ 
\frac{(i_1+i_2-j)(i_2-i_1-j)}{4k_1}}.$$
\end{cor}

\newcounter{a}
\setcounter{a}{2}


\begin{thebibliography}{99}

\bibitem[B]{B}
J. Bagger, D. Nemeschansky, S. Yankielowicz, 
{\it Virasoro algebras with central charge $c>1$},
Phys. Rev. Lett. {\bf 60} (1988), no. 5, 389-392.

\bibitem[BGG]{BGG}
I.N. Bernstein, I.M. Gel'fand, S.I. Gel'fand, {\it Differential operators on
the base affine space and a study of $\g$-modules}, in: Lie groups and
their representations (I.M. Gelfand ed.), Summer school of the Bolyai Janos
Math. Soc., Halsted Press, 1975, 21-64.

\bibitem[DJKMO]{DJKMO}
E. Date, M. Jimbo, A. Kuniba, T. Miwa, M. Okado, 
{\it Exactly solvable SOS models: local heights probabilities and theta 
function identities}, Nucl. Phys. B {\bf 290} (1987), no. 2, 231-273.

E. Date, M. Jimbo, A. Kuniba, T. Miwa, M. Okado, 
{\it Exactly solvable SOS models \Roman{a}: proof of star-triangle relation
and combinatorial identities}, Adv. Stud. in Pure Math. 
{\bf 16} (1988), 17-122.

\bibitem[FF1]{FF1}
B.Feigin, E.Feigin,
{\it Q-characters of the tensor products in $\slt$ case,}
Mosc. Math. J. {\bf 2}, no. 3, 567-588.

\bibitem[FF2]{FF2}
B.Feigin, E.Feigin,
{\it Integrable $\slth$-modules as infinite tensor products.},
Fundamental mathematics today (S.Lando and O.Sheinman eds.),
Independent University of Moscow, 2003, pp. 304-334, (in honour
of the $10$th anniversary of the Independent University of Moscow).

\bibitem[FF3]{FF3}
B. Feigin, E. Feigin, {\it Schubert varieties and the fusion products},
Publ. Res. Inst. Math. Sci., Kyoto Univ. {\bf 40} (2004), no. 3, 625-668.

\bibitem[FF4]{FF4}
B. Feigin, E. Feigin, {\it Homological realization of restricted Kostka
polynomials}, Int. Math. Res. Not. 2005, no. 33, pp. 1997-2029.

\bibitem[FFu]{FFu}
B. Feigin, D. Fuchs, {\it Verma modules over a Virasoro algebra},  
Funct. Anal. Appl. {\bf 17} (1983), no. 3, 91-92.

\bibitem[FJKLM]{FJKLM}
B.Feigin, M.Jimbo, R.Kedem, S.Loktev, T.Miwa, {\it Spaces of coinvariants
and fusion product \Roman{a}. $\slth$ character formulas in terms of Kostka
polynomials},
J. Algebra {\bf 279} (2004), no. 1, 147-179.


\bibitem[FL]{FL}
B.Feigin, S.Loktev,
{\it On generalized Kostka polynomials and quantum Verlinde rule,}
Differential topology, infinite-dimensional Lie algebras and applications,
Amer. Math. Soc. Transl. Ser.  2, vol. 194, American Mathematical Society,
Rhode Island, 1999, 61-79.



\bibitem[GKO]{GKO}
P. Goddard, A. Kent, D. Olive, {\it Unitary representations of the Virasoro
and super-Virasoro algebras}, Comm. Math. Phys., {\bf 103} (1986), no.1,
105-119.

\bibitem[GL]{GL}
H. Garland, J. Lepowsky, {\it Lie algebra homology and the Macdonald-Kac
formulas}, Invent. Math. {\bf 34} (1976), 37-76.


\bibitem[Kac]{Kac}
V. Kac, {\it Infinite dimensional Lie algebras}, 3rd ed.,
Cambridge University Press, Cambridge, 1990.

\bibitem[K]{K}
D. Kastor, E. Martinec, Z. Qiu, {\it Current algebra and conformal decrete
series}, Phys. Lett. B{\bf 200} (1988), no. 4, 434-440. 

\bibitem[R]{R} F. Ravanini, {\it An infinite class of new conformal field
theories with extended algebras},  Mod. Phys. Lett. A{\bf 3} (1988), no. 4, 
397-412.

\bibitem[RC]{RC}
A. Rocha-Caridi, {\it Vacuum vector representations of the Virasoro algebra},
Vertex Operators in Mathematical Physics (Berkeley, Calif., 1983),
Math. Sci. Res. Inst.  Publ., vol. 3, Springer, New York, 1985, pp. 451-473.


\bibitem[S1]{S1}
A. Schilling, {\it Multinomials and polynomial bosonic forms for the branching
functions of the
$\widehat{\mathfrak{su}}(2)_M\times \widehat{\mathfrak{su}}(2)_N/
\widehat{\mathfrak{su}}(2)_{M+N}$ conformal coset models}, 
Nucl. Phys. B{\bf 467} (1996), 247-271.    


\bibitem[S2]{S2}
A. Schilling, {\it Polynomial fermionic forms for the branching functions
of the rational coset conformal field theories
$\widehat{\mathfrak{su}}(2)_M\times \widehat{\mathfrak{su}}(2)_N/
\widehat{\mathfrak{su}}(2)_{M+N}$}, Nucl. Phys. B{\bf 459} (1996), 393-436.

\bibitem[SS]{SS}
A. Schilling, M. Shimozono, {\it Fermionic formulas for level-restricted
generalized Kostka polynomials and coset branching functions}, Comm.
Math. Phys. {\bf 220} (2001), no. 1, 105-164.



\end{thebibliography}
\end{document}